\documentclass{article}
\usepackage{amsmath}
\font\tengoth=eufm10
\font\sevengoth=eufm7
\font\fivegoth=eufm5
\newfam\gothfam
\textfont\gothfam=\tengoth \scriptfont\gothfam=\sevengoth
\scriptscriptfont\gothfam=\fivegoth

\newtheorem{theorem}{Theorem}[section]

\newtheorem{corollary}[theorem]{Corollary}

\def\blacksquare{\hbox to .60em{\vrule width .60em height .60em}}
  \font\bb=msbm10 
\catcode`\@=12  
\def\é{\'e}
\def\{\`e}
\def\?{\`a}
\def\{\`u}
\def\{\c c}
\def\hb {\hfil \break}
\def\n {\vskip 0.2cm \noindent }
\def\scirc{\,{\raise 0.8pt\hbox{$\scriptstyle\circ$}}\,}
\def\ins{\,{\raise 0.2cm \hbox{ $\scriptstyle \circ$}}\,}

\def  \é{\'e}
\def\è{\`e}
\def\à{\`a}
\def\ù{\`u}
\def\ç{\c c$\!\!\!$}
 
 \date{ }
 \title{Ordinary algebraic curves and related topics}
\author{   Youssef Hantout, and Daniel Lehmann}

\begin{document}
  \centerline {\bf Courbes ordinaires de genre maximal}
 %\centerline {\bf  Rank and Blaschke curvature}
\medskip

  \centerline  {Youssef Hantout et Daniel Lehmann.}
  
  \bigskip
  
  \rightline{\à la m\émoire de L. Gruson}
  
  % \rightline {\it
 % Work in progress ; not to be diffused.}
  
   \bigskip

\section{Introduction}

Les courbes dont il s'agit dans cet article sont   des courbes alg\ébriques dans l'espace projectif complexe {\bb P}$_n$ ($n\geq 3)$, toujours  suppos\ées irr\éductibles, et propres (c'est-\à-dire non incluses dans un sous-espace projectif strict de {\bb P}$_n$).

Les courbes alg\ébriques \emph{ordinaires}\footnote{Ce concept a \ét\é sugg\ér\é \à l'origine  par celui de  {\it tissu ordinaire} ([CL]). } de degr\é $d$ sont   les courbes $\Gamma$ pour lesquelles les $d$ points d'une section hyperplane g\én\érique $H\cap \Gamma$ sont ``en position g\én\érale''. Pour $n=3$ par exemple, cela signifie que 3 points de  $H\cap \Gamma$ ne sont jamais align\és, 6 points ne sont jamais sur une m\^eme conique, etc....\hb Plus g\én\éralement,  pour $n$ quelconque,
notons $$c(n,h):=(n-1+h)!/(n-1)!h!$$   la dimension   de l'espace {\bb C}$_h[H]$ des polyn\^omes homog\ènes de degr\é $h$ en $n$ variables.
Dire que la courbe est ordinaire  signifie alors que, pour tout entier $h$ ($h\geq 1)$, $c(n,h)$ points distincts   d'une section hyperplane g\én\érique  $H\cap \Gamma$ n'appartiennent pas \à une m\^eme hypersurface alg\ébrique de degr\é $h$ dans l'hyperplan  $H$. 
%[Il revient au m\^eme d'\énoncer que l'application naturelle de restriction 
%$$H^0\bigl(H,{\cal O}_H(h)\bigr)\to H^0\bigl(H\cap \Gamma,{\cal O}_{H\cap \Gamma}(h)\bigr)$$
%est de rang maximal, et ceci  quel que soit l'entier $h\geq 1$]

Les courbes alg\ébriques \emph{arithm\étiquement de Cohen-Maccaulay} (acm en abr\ég\é)  sont   les courbes $\Gamma$ pour lesquelles 
les seules hypersurfaces $S$ d'un hyperplan g\én\érique $H$ qui contiennent $H\cap \Gamma$ sont les intersections $H\cap S'$ avec $H$ d'une hypersurface $S'$ de 
{\bb P}$_n$ qui contient  $\Gamma$.

Soit  $d$ le degr\é d'une courbe alg\ébrique $\Gamma$ de {\bb P}$_n$ ($d\geq n$), $k_0$ l'entier ($\geq 1$) tel que 
$$c(n,k_0)\leq d <c(n,k_0+1),$$ et 
$$\pi'(n,d):=\sum_{h=1}^{k_0} \bigl(d-c(n,h)\bigr)\ \ \ \bigl(=k_0.d-c(n+1,k_0)+1\bigr).$$
Avec L. Gruson, nous avons d\émontr\é le 
 
\n {\bf Th\éor\ème 1 }([GHL]) : 

\n {\it L'entier $\pi'(n,d)$ est \à la fois 

- une borne sup\érieure du genre arithm\étique des courbes ordinaires de degr\é $d$ dans {\bb P}$_n$.

- une borne inf\érieure du genre arithm\étique des courbes acm de degr\é $d$ dans {\bb P}$_n$,

\n En outre, les courbes ordinaires  de genre maximal   sont acm, et les courbes acm de genre minimal   sont  ordinaires.}

\n {\it D\émonstration } (non publi\ée dans [GHL]) : 
%{\it D\émonstration :} 

Pour tout ensemble   $A=\{m_1,\cdots,m_d\}$ de $d$  points  distincts dans un hyperplan projectif   $H$ de {\bb P}$_n$, les polyn\^omes homog\ènes de degr\é $h$ sur $H$ qui s'annulent sur $A$ sont solutions d'un syst\ème lin\éaire homog\ène  $\Phi_h(A)$ de $d$ \équations \à $c(n,h)$ inconnues. Le rang  $\lambda_h(A)$ de ce syst\ème est donc  au plus \égal \à l'entier 
$$ 
  \lambda_0(d,n,h) = \hbox{ min }\bigl(d, c(n,h)\bigr). 
  $$
   %[On observe que $\sum_{h=1}^\infty \bigl(d-\lambda_0(d,n,h)\bigr)=\pi'(n,d)$].
\pagebreak

  D'autre part (cf. par exemple [HE]), le genre arithm\étique $g(\Gamma)$ d'une courbe alg\ébrique  $\Gamma$  dans {\bb P}$_n$, propre et de degr\é $d$,  
   v\érifie  l'in\égalit\é 
  $$g(\Gamma)\leq \sum_{h=1}^\infty \bigl(d-\lambda_h(H\cap \Gamma)\bigr)$$
  pour toute section hyperplane g\én\érique $H$. En outre, cette in\égalit\é  devient une \égalit\é ssi $\Gamma$ est acm. 
  
 % \pagebreak
   
 Puisque $\lambda_h(H\cap \Gamma)= \lambda_0(d,n,h)$ pour une courbe ordinaire, le genre arithm\étique d'une telle  courbe  est  major\é par $\sum_{h=1}^\infty \bigl(d-\lambda_0(d,n,h)\bigr)$, c'est-à-dire par $\pi'(n,d)$. Si cette borne est atteinte, c'est que l'in\égalit\é est en fait une \égalit\é: la courbe est donc  acm. 
 
  Puisque $d-\lambda_h(H\cap \Gamma)\geq d- \lambda_0(d,n,h)$, 
 le genre arithm\étique d'une courbe acm est  minor\é par $\pi'(n,d)$. Si la borne est atteinte, c'est que $\lambda_h(H\cap \Gamma)=\lambda_0(d,n,h)$ pour tout $h$ : la courbe est donc ordinaire.

\rightline{QED}
Dans [GHL], nous avons montr\é que, pour $n=3$, il existait en tout degr\é $d$ ($d\geq 3$) des courbes \à la fois acm et otrdinaires, donc de genre arithm\étique $\pi'(3,d)$.  [On pouvait  m\^eme les choisir lisses, et leur ensemble constitue  une composante irr\éductible du sch\éma de Hilbert ${\cal H}_{\pi'(3,d),d}$].

Nous  nous proposons,  dans  cet article, de montrer l'existence de courbes \à la fois acm et ordinaires,  $\bigl($donc de genre arithm\étique $\pi'(n,d)\bigr)$, 

- pour tout $n$ ($n\geq 3$) et tout $d$ calibr\é (c'est-\à-dire \égal \à $c(n,k_0)$, $k_0\geq 1$),

- et plus g\én\éralement \à condition de savoir construire   une suite exacte de fibr\és   d'une certaine forme que nous allons pr\éciser  (ce qui est en particulier toujours possible  quel que soit $d$, si $n=3$). 

\section{Conditions suffisantes pour qu'une suite exacte de fibr\és soit la r\ésolution d'une courbe acm et ordinaire }

Supposons disposer d'une suite exacte de fibr\és vectoriels holomorphes au dessus de    {\bb P}$_n$ (not\é parfois {\bb P} en abr\ég\é dans la suite).
$$(*)\hskip 1cm 0\to E_{N}\buildrel{\partial _{N-1}}\over\longrightarrow  E_{N-1}\to \cdots\to E_{i+1}\buildrel{\partial _i}\over\longrightarrow E_{i}\to \cdots \to E_1\buildrel{\partial _0}\over\longrightarrow E_0= {\cal O}_{\hbox {\bb P} }  $$
o\ù chaque fibr\é $E_i$ est de la forme    $\oplus _{j=1}^{J_i}{\cal O}_{\hbox {\bb P} }(-a_{j,i})$ pour une famille d'entiers positifs  $a_{1,i},\cdots, a_{J_i,i}$. On suppose   $a_{j,i}\geq a_{j',i'}$ pour   $i>i'$, de fa\c con que tous les   morphismes de cette suite  exacte puissent \^etre repr\ésent\és par des   matrices dont les coefficients non nuls\footnote{Si le degr\é de ces coefficients non nuls est strictement positif, on dit alors que la suite $(*)$ est \emph{minimale}. C'est en particulier ce qui se produit si l'on a des in\égalit\és strictes ($a_{j,i}> a_{j',i'}$ pour   $i>i'$),  la r\éciproque \étant  fausse.} sont des polyn\^omes homog\ènes    par rapport aux   coordon\ées homog\ènes  $(X_0,\cdots,X_n)$ de  {\bb P}$_n$. 

  \n Le conoyau   de $\partial _0$   est l'anneau  ${\cal O}_{ \Gamma }$ d'un sous-ensemble alg\ébrique   $\Gamma$ de  {\bb P}$_n$.

\n $(i)$  Pour que $\Gamma$ soit une courbe de degr\é $d$ et genre arithm\étique $g$,  il faut et il suffit que sa fonction de     Hilbert  
$$H_{\Gamma}(h)=c(n+1,h)-\sum_{i=1}^{N}(-1)^i\sum_{j=1}^{J_i}c(n+1,h-a_{j,i})$$ co\"\i ncide     avec le polyn\^ome de    Hilbert      $P(h)= hd-g+1$ pour  $h$   suffisamment grand.

\n $(ii)$ Et pour que cette courbe soit acm,  il suffit que la longueur  de la suite soit \égale \à la codimension $n-1$, soit 
  $$N=n-1.$$
  Posant  :
$$S_p=\sum_{i=1}^{n-1}(-1)^i\sum_{j=1}^{J_i} (a_{j,i})^p $$
pour tout entier   $p\geq 0$, la condition $(i)$ s'\écrit : 
$$S_0=1\hbox{ (i.e. } \sum_{i=1}^{n-1}(-1)^i J_i=1), \hskip 1cm S_p=0 \hbox{ pour  } 1\leq p\leq n-2,\hskip 1cm S_{n-1}=(-1)^n (n-1)!\ d,$$
et  
$$S_n=(-1)^n n!\Bigl(g-1+\frac{n+1}{2}\ d\Bigr).$$

\medskip

\n  Cherchons s'il est possible que  les entiers   $a_{j,i}$ soient    r\épartis parmi    $n$ entiers  $b_1,\cdots,b_n$ tels que  $$0<b_1<b_2<\cdots<b_n.$$
Si l'on   impose   des  {\bf   in\égalit\és strictes}  $a_{j,i}>a_{j',i'}$ pour  $i>i'$, on peut proc\éder   a priori de $n-1$ fa\c cons distinctes   : il  existe un  indice  $i_0$ ($1\leq i_0\leq n-1$) tel que \hb  
$$a_{j,i}=b_{i},\hbox{ pour  } i<i_0\hskip .5cm  a_{j,i}=b_{i+1}\hbox{  pour  } i>i_0 ,\hskip .5cm  \hbox{ et   } a_{j,i_0}=b_{i_0}  \hbox{ ou } b_{i_0+1}.$$ 

\n Notons  $x_i$ le nombre de   $a_{p,q}$  \égaux \à  $b_i$ ($1\leq i\leq n$), de sorte que
$J_i=x_i$  si $i<i_0$, $J_i=x_{i+1}$ si $i>i_0$, et $J_{i_0}=x_{i_0}+x_{i_0+1}$. Chacune de ces r\épartitions conduit au   syst\ème de  $n$ \équations  lin\éaires 
$$S_0=1 , \hskip 1cm S_p=0 \hbox{ pour   } 1\leq p\leq n-2,\hskip 1cm S_{n-1}=(-1)^n (n-1)!\ d $$
\à   $n$ inconnues  $x_i$ que l'on peut \écrire sous la forme 

$$\begin{pmatrix}
1&1&\cdots&1& \cr
b_1&b_2&\cdots&b_n& \cr 
(b_1)^2&(b_2)^2&\cdots&(b_n)^2& \cr 
\vdots &\vdots& &\vdots&\cr
(b_1)^{n-2}&(b_2)^{n-2}&\cdots&(b_n)^{n-2}& \cr 
(b_1)^{n-1}&(b_2)^{n-1}&\cdots&(b_n)^{n-1}& \cr 
 
\end{pmatrix}.\begin{pmatrix}
\overline x_1 \cr
\overline x_2 \cr 
\vdots \cr 
\overline x_{n-2}\cr
\overline x_{n-1}\cr 
\overline x_{n} \cr 
 
\end{pmatrix}=\begin{pmatrix}
1 \cr
0\cr 
\vdots \cr 
0\cr
0\cr 
 (-1)^n(n-1)!\ d \cr 
 
\end{pmatrix} ,$$
\à condition de poser  $\overline x_i=-x_i$ pour $1\leq i\leq i_0$, et $\overline x_i=x_i$ pour $i_0+1\leq i\leq n$.

\n  Ces syst\èmes sont  crameriens    (matrice  de Vandermonde relative aux nombres   $b_i$), et ont donc tous une   solution unique.  

\n Mais nous avons aussi besoin de choisir l'indice   $i_0$ de fa\c con que tous les nombres $x_i$ de la solution soient des entiers, et des entiers positifs, strictement    pour  $i<i_0$, $i>i_0+1$, et pour au moins l'un des  indices  $i=i_0$ ou $i_0+1$. Nous allons voir que c'est 
toujours possible, et d'une seule fa\c con. Mais admettons le  provisoirement. 
 On obtient alors une courbe  acm   $\Gamma$ de degr\é  $d$, et dont le genre  $g$ peut se calculer en \écrivant que la  $(n+1)^{i\grave{e}me}$  \équation $S_n= (-1)^n n!\bigl(g-1+\frac{n+1}{2} d \bigr)$ est compatible avec les pr\éc\édentes, soit : 
 $$\begin{vmatrix}
1&1&\cdots&1&&1 \cr
b_1&b_2&\cdots&b_n&&0 \cr 
(b_1)^2&(b_2)^2&\cdots&(b_n)^2&&0 \cr 
\vdots &\vdots& &\vdots&&\vdots\cr
(b_1)^{n-2}&(b_2)^{n-2}&\cdots&(b_n)^{n-2}&&0 \cr 
(b_1)^{n-1}&(b_2)^{n-1}&\cdots&(b_n)^{n-1}&&(-1)^n(n-1)!\ d \cr 
(b_1)^{n }&(b_2)^{n }&\cdots&(b_n)^{n }&&(-1)^n n!\bigl(g-1+\frac{n+1}{2} d \bigr)\cr 
\end{vmatrix}=0\ .$$
 D\éveloppant le d\éterminant ci-dessus  par rapport \à la derni\ère colonne, et notant  $D$ le d\éterminant de   Vandermonde construit avec les nombres   $(b_i)_{1\leq i\leq n}$,  la derni\ère  \équation devient   :
  $$ D.\biggl( \sigma_n  -  \sigma_1     (n-1)!\   d +       n!\  \Bigl(g-1+\frac{n+1}{2} d\Bigr)\biggr)=0,$$
  dans laquelle  $\sigma_1:=\sum_{i=1}^n b_i  $,  et $\sigma_n:=\prod_{i=1}^n b_i  $, d'o\ù :
  $$g=\Bigl(\frac{\sigma_1}{n}-\frac{n+1}{2}\Bigr)\ d-\frac{1}{n!}\ \sigma_n+1.$$
  Rappelons l'\égalit\é $$\pi'(n,d)=k_0.d-\frac{1}{n!}\prod_{i=1}^n (k_0+i)+1.$$  Il suffit donc que 
  $\frac{\sigma_1}{n}-\frac{n+1}{2}=k_0$ et $\sigma_n=\prod_{i=1}^n (k_0+i)$ pour que $g$ soit \égal \à $\pi'(n,d)$ : d'apr\ès le th\éor\ème 1 rappel\é dans l'introduction,  la courbe acm sera   de genre minimal, donc ordinaire. 
  Ces conditions sont visiblement r\éalis\ées si l'on d\éfinit 
  $$b_i:=k_0+i \hbox { pour tout }i=1, \cdots,n.$$ 
  Montrons alors que l'un des syst\èmes cram\ériens convient pour ces valeurs des $b_i$, \à condition de   choisir l'indice $i_0$   en fonction de la position de $d$ dans l'intervalle $\bigl[ c(n,k_0), c(n,k_0+1)\bigr [$.
On d\éfinit pour cela la suite d\écroissante $$d_1=c(n,k_0+1)>d_{n-1}>\cdots>d_n=c(n,k_0)$$ en posant :
$$d_i=\frac{1}{(n-1)!}\prod_{ 1\leq j\leq n,\  j\neq i }(k_0+j).$$
  On trouve alors :
 $$\overline x_i=\frac{(n-1)!}{(i-1)!\ (n-i)!}\ (d_i-d),$$
% o\ù l'on a pos\é: 
 %$$\Delta_i=\prod_ {1\leq j\leq n,  (j\neq i)} |i-j| \  \ \ \bigl(=(i-1)!\ (n-i)!\bigr).$$ 

\n {\bf Lemme  :}

\n {\it Les nombres   $x_i$ sont tous des entiers, et des entiers positifs si l'on choisit pour  $i_0$ $(1\leq i_0\leq n-1)$ l'indice tel que $d\in [d_{i_0+1},d_{i_0}[$.   

}

\n{\it D\émonstration :} 
Le nombre  $\overline x_i$  est en effet \égal au nombre entier   $$c(i,k_0).c(n-i-1,k_0+i)-\frac{(n-1)!}{(i-1)!(n-i)!}\ d\  ,$$ et $x_i\geq 0$, grace au choix de   $i_0$. 
\rightline{QED}

Nous avons donc d\émontr\é le 

\n {\bf Th\éor\ème 2 :}

{\it 
S'il existe des morphismes ${\partial _i}$ tels que la suite de fibr\és vectoriels 
$$(*)\hskip 1cm 0\to E_{n-1}\buildrel{\partial _{n-2}}\over\longrightarrow  E_{n-2}\to \cdots\to E_{i+1}\buildrel{\partial _i}\over\longrightarrow E_{i}\to \cdots \to E_1\buildrel{\partial _0}\over\longrightarrow E_0= {\cal O}_{\hbox {\bb P} }  $$
soit exacte, o\ù l'on a pos\é : 
$$\begin{matrix}E_i =&\hskip -3cm  \frac{(n-1)!}{\Delta_{i}}\ (d_{i}-d)\ {\cal O}_{\hbox{\bb P}}\bigl(-(k_0+i)\bigr)   \hbox {\hskip 1cm pour  }  1\leq i< i_0\ ,\cr
&  \cr
E_{i_0} =& \frac{(n-1)!}{\Delta_{i_0}}\ (d_{i_0}-d)\ {\cal O}_{\hbox{\bb P}}\bigl(-(k_0+i_0)\bigr)\oplus \frac{(n-1)!}{\Delta_{i_0+1}}\ (d-d_{i_0+1})\ {\cal O}_{\hbox{\bb P}}\bigl(-(k_0+i_0+1)\bigr)  \ ,\cr
&   \cr
E_i =& \hskip -2cm\frac{(n-1)!}{\Delta_{i+1}}\ (d-d_{i+1})\ {\cal O}_{\hbox{\bb P}}\bigl(-(k_0+i+1)\bigr)  \hbox { \hskip 1cm pour } i_0< i\leq n-1,\cr
\end{matrix}, $$ 
le conoyau de $\partial_0$ est alors l'anneau d'une courbe \à la fois acm et ordinaire  $($donc de genre arithm\étique $\pi'(n,d))$.}
%where $N  {\cal O}_{\hbox{\bb P}} (-j )$ states for ${\cal O}^{\oplus N}_{\hbox{\bb P}} (-j )$ (the Whitney sum of $N$ copies of ${\cal O} _{\hbox{\bb P}} (-j )$). It is not hard to prove that the bundles above are convenient. 
 
  \n {\bf Remarque :} Si l'on permet des  {\bf  in\égalit\és  larges} $a_{j,i}\geq a_{j',i'}$  pour    $i>i'$, ce qui suffirait pour que l'on puisse esp\érer d\éfinir des morphismes $\partial_i$,  la r\épartition possible    des   copies de  ${\cal O} _{\hbox {\bb P} }(-(k_0+j))$ n'est alors plus   unique.  
  
  On peut \évidemment  cr\éer des "redondances" en ajoutant un m\^eme fibr\é
   $u {\cal O} _{\hbox {\bb P}}  (-(k_0+j))$   \à deux   fibr\és cons\écutifs   $E_{i+1}$ et  $E_{i}$ d'une r\ésolution donn\ée,  et en modifiant $\partial_i$ \à l'aide d'un isomorphisme entre les deux copies de  $u. {\cal O} _{\hbox {\bb P}}  (-(k_0+j))$ : ce faisant,  en effet, on ne changera pas les nombres   $S_p$. La nouvelle r\ésolution     ne sera certainement pas minimale. 
  
Mais il se peut  aussi que    des r\ésolutions minimales  puissent \^etre obtenues  avec  in\égalit\és larges (voir le deuxi\ème exemple donn\é pour $n=4$).

   \section{Exemples}

\subsection{Courbes  calibr\ées}

Nous dirons qu'une courbe alg\ébrique  $\Gamma$ de degr\é $d$  dans {\bb P}$_n$  est \emph{calibr\ée},  s'il  
 existe un entier $k_0$   au moins \égal \à 1 tel que $d=c(n,k_0) $.   Par exemple, pour $n=3$ (resp. 4), les courbes calibr\ées sont les courbes gauches de degr\é 3,  6, 10, 15, 21 .....(resp. 4, 10, 20, 35, 56 .....).  
 
%On a alors , $d=d_n$, $i_0=n-1$, $E_i=F_i$ quel que soit $i=1,\cdots,n-1$, et $F_n=0$.  Le complexe d'Eagon-Northcott d\éfini dans la section 2 ci-dessus fournit la suite exacte (*) : il suffit de  v\érifier   que les nombres $x_i$ d\éfinis dans le cas g\én\éral co\"\i ncident avec ceux d\éfinis dans la section 2.

%Notons  $F_i$ le fibr\é vectoriel   $x_i. {\cal O}_{\hbox{\bb P}}\bigl(-(k_0+i)\bigr)$ de base  {\bb P}={\bb P}$_n$,
%\égal \à  la somme de  Whitney ${\cal O}^{\oplus |x|}_{\hbox{\bb P}} (-j )$    de $|x|$ copies de la puissance tensorielle j-i\ème ${\cal O} _{\hbox{\bb P}} %(-j )$ du fibr\é tautologique ${\cal O} _{\hbox{\bb P}} (-1 )$ en droites complexes. 

\n A partir  d'un morphisme de fibr\és vectoriels\footnote{Il  revient au m\^eme de se donner une matrice de taille $(k_0+1)\times (k_0 +n-1)$ dont tous les coefficients sont des combinaisons lin\éaires des coordonn\ées homog\ènes de {\bb P}$_n$.}
 $$\varphi:(k_0+n-1).{\cal O}_{\hbox{\bb P}} (-1)\  \to\ (k_0+1).{\cal O}_{\hbox{\bb P}} (0) ,$$
on   d\éfinit   le complexe d'Eagon-Northcott 
 $$0\to\bigwedge^{k_0+n-1}V\otimes S^{n-2}W^*\to \cdots\to \bigwedge^{k_0+i+1}V\otimes S^{i}W^*\buildrel{\partial_i}\over\longrightarrow \bigwedge^{k_0+i}V\otimes S^{i-1}W^*  \cdots \to  \bigwedge^{k_0+1}V\buildrel{\partial_0}\over\longrightarrow {\cal O} _{\hbox{\bb P}} $$ 
 \n  o\ù $V$ $($resp. $W)$ d\ésigne le fibr\é vectoriel $\bigl(k_0+(n-1)\bigr).{\cal O}_{\hbox{\bb P}}(-1)$ $\bigl($resp.  $(k_0+1).{\cal O}_{\hbox{\bb P}}(0)\bigr)$, 
 
 \n  $\bigwedge$ et $S$ d\ésignent respectivement des puissances ext\érieures et sym\étriques., 
 
\n    $\partial_0$ est \égal au morphisme  $\wedge^{k_0+1}\varphi$ de $\bigwedge^{k_0+1}V$ dans $\bigwedge^{k_0+1}W\cong {\cal O} _{\hbox{\bb P}}$.
 
 \n  pour $i>0$, on d\éfinit $\partial_i:\bigwedge^{k_0+i+1}V\otimes S^{i}W^*\to \bigwedge^{k_0+i}V\otimes S^{i-1}W^* $ par la formule 
 $$\partial_i\bigl(v_{1}\wedge  \cdots\wedge v_{k}\wedge\cdots\wedge v_{k_0+i+1}\otimes \epsilon_1. \epsilon_2.\cdots\epsilon_u.\cdots\epsilon_i\bigr)$$ 
  $$= \sum_k\sum_u (-1)^k  <\epsilon_{u}\scirc \varphi,v_k>.(v_{1}\wedge  \cdots\widehat{v_{k}}\cdots\wedge v_{k_0+i+1}\otimes \epsilon_1. \epsilon_2.\cdots\widehat{\epsilon_u}.\cdots\epsilon_i),$$
 dans laquelle 
 
 $v_{1}\wedge  \cdots\wedge v_{k}\wedge\cdots\wedge v_{k_0+i+1}\in \bigwedge^{k_0+i+1}V$  d\ésigne un produit ext\érieur  d'\él\éments $v_k$ dans une fibre de $V$ au dessus d'un point de {\bb P}$_n$,  
 
  $ \epsilon_1. \epsilon_2.\cdots\epsilon_u.\cdots\epsilon_i$ d\ésigne un produit sym\étrique de formes lin\éaires $\epsilon_u$ (pas n\écessairement toutes distinctes) sur  la fibre de $ W$ au dessus  du m\^eme point de {\bb P}$_n$,
  
 tandis que $\widehat{v_k}$ (resp. $\widehat{\epsilon_u}$)  signifient que $v_k$ (resp. $\epsilon_u$) ont \ét\é omis. 
  
\n  Le morphisme $\partial_0$  a   pour image l'id\éal $I(\varphi)$ de ${\cal O} _{\hbox{\bb P}}$ engendr\é par les mineurs de taille maximale $(k_0+1)\times (k_0+1)$ de la matrice $\varphi$. On dira que $\varphi$ est \emph{g\én\érique} s'il
 existe une suite r\éguli\ère de longueur $n-1$ dans cet id\éal. De tels morphismes $\varphi$ g\én\ériques existent toujours (cf. [E]).

 \n On obtient toujours un complexe ($\partial_{i-1}\scirc\partial_i=0$) et, si $\varphi$ est g\én\érique,  ce complexe 
est  une suite exacte de fibr\és vectoriels (cf. par exemple [E]).

 \n {\bf Th\éor\ème 3 :}
 
\n  {\it Si  $\varphi$ est {g\én\érique}, le sous-ensemble alg\ébrique $\Gamma$ de {\bb P}$_n$ d\éfini par l'id\éal $I(\varphi)$ est une courbe  de degr\é $d=c(n,k_0)$, \à la fois ordinaire et acm, et par cons\équent de genre arithm\étique  $g(\Gamma)=\pi'(n,d)$.\hb 
Il existe donc de telles courbes quels que soient $n$ $(n\geq 3)$ et $k_0$ $(k_0\geq 1)$.}
 
\n {\it D\émonstration :}

%On a alors , $d=d_n$, $i_0=n-1$, $E_i=F_i$ quel que soit $i=1,\cdots,n-1$, et $F_n=0$.  Le complexe d'Eagon-Northcott d\éfini dans la section 2 ci-dessus fournit la suite exacte (*) : il suffit de  v\érifier   que les nombres $x_i$ d\éfinis dans le cas g\én\éral co\"\i ncident avec ceux d\éfinis dans la section 2. 

\n Remarquons tout d'abord 
l'existence, pour tout $i=1,\cdots,n-1$, d'un isomorphisme de  fibr\és vectoriels 
 $$\bigwedge^{k_0+i}V\otimes S^{i-1}W^*\cong y_i.   {\cal O}_{\hbox{\bb P}} (-(k_0+i) ), \hbox{ o\ù  }y_{i}=\begin{pmatrix}k_0+n-1\cr k_0+i\end{pmatrix}.\begin{pmatrix}k_0+i-1\cr i-1\end{pmatrix}.$$
En effet,  $\bigwedge^{k_0+i}V$ est isomorphe \à $\begin{pmatrix}k_0+n-1\cr k_0+i\end{pmatrix}.{\cal O}_{\hbox{\bb P}}\bigl(-(k_0+i)\bigr)$, 
 et $ S^{i-1}W^*$  \à $c(k_0+1,i-1). {\cal O}_{\hbox{\bb P}}(0)$.
 
\n Puisque  $d=d_n$ et  $i_0=n-1$,  il suffit donc de   v\érifier   que les nombres $x_i$ et $y_i$ sont \égaux pour prouver que le complexe d'Eagon-Northcott est du type d\écrit dans le th\éor\ème 2. On v\érifie ais\ément que chacun de ces deux nombres est \égal \à 
$$\frac{1}{(k_0+i)}.\frac{(k_0+n-1)!}{k_0!(i-1)!(n-i-1)!} .$$

\rightline{QED}

%\subsection {\bf Cas o\ù $d_{i_0+2}$ est un  entier, et  $d=d_{i_0+2}$ :}
 
\subsection {\bf Cas  $n=3$ :}  C'est le cas \étudi\é dans [GHL]. 
On a ici  :
 
 \n  $d_1=\frac{   1}{2}(k_0+2)(k_0+3) \ \bigl(=c(3,k_0+1)\bigr)$,\hb 
$d_2=\frac{ 1 }{2} (k_0+1)(k_0+3) $, \hb
  $d_3=\frac{   1}{2}(k_0+1)(k_0+2) \ \bigl(=c(3,k_0)\bigr)$.
  
 \n  Posant $r=d-c(3,k_0)$ et $t=c(3,k_0+1)-d$  $($de sorte que   $r+t=k_0+2)$, et $d\in [d_3,d_2] $ ou $d\in [d_2,d_1] $ selon que $t\geq r+1$ ou 
 $t\leq r+1$. 
 
 \n Quel que soit $d$ $(d\geq 3)$, il existe des  suites exactes du type voulu, ainsi qu'il r\ésulte de  [GHL] :    
 
\n {\bf Th\éor\ème 4} ([GHL])

\n   {\it $(i)$ Pour tout entier   $d$ $(d\geq 3)$, la famille des courbes   de degr\é $d$  dans  $ \hbox {\bb P}_3$ qui sont \à la fois acm et ordinaires est non vide,  et  constitue une   composante irr\éductible  du sch\éma de   Hilbert ${\cal H}_{d,\pi'(3,d)}$. 
 
 \n $(ii)$    Soit $M$   un   morphisme injectif de fibr\és vectoriels   $($dont tous les  coefficients non nuls  sont des    polyn\^omes homog\ènes de  degr\é  1 ou 2 $)$   : 

 - si \ $t\geq r+1$,
 $$r{\cal O}_{\hbox {\bb P}  }\bigl(-(k_0+3)\bigr)  \oplus (t-r-1 ){\cal O}_{\hbox {\bb P}  }\bigl(-(k_0+2)\bigr)                  \buildrel M \over
  \longrightarrow t{\cal O}_{\hbox {\bb P}   }\bigl(-(k_0+1)\bigr) $$       
   
 - si $t\leq r+1$, 
$$r{\cal O}_{\hbox {\bb P} _3}\bigl(-(k_0+3)\bigr)  \buildrel M \over \longrightarrow   (r+1-t ){\cal O}_{\hbox {\bb P}_3 }\bigl(-(k_0+2)\bigr)\oplus  t{\cal O}_{\hbox {\bb P}_3 }\bigl(-(k_0+1)\bigr)  $$   

- et plus g\én\éralement :
 $$ r{\cal O}_{\hbox {\bb P}  }\bigl(-(k_0+3)\bigr)  \oplus (t-r-1+u){\cal O}_{\hbox {\bb P}}\bigl(-(k_0+2)\bigr)                \buildrel M \over \longrightarrow   u {\cal O}_{\hbox {\bb P}}\bigl(-(k_0+2)\bigr)\oplus t{\cal O}_{\hbox {\bb P}}\bigl(-(k_0+1)\bigr) $$        
 o\ù $u\geq 0$ d\ésigne un entier positif tel que $t-r-1+u\geq 0$. 
 
    Les cofacteurs de taille maximale de     $M$ engendrent alors l'id\éal  d'une courbe ordinaire et acm de degr\é  $d$. 
 
  \n $(iii)$ Il est en outre possible de choisir   $M$ de fa\c con que la courbe soit lisse. }
  
 \n La  partie $(ii)$ est un cas particulier du th\éor\ème 2.  Nous renvoyons \à [GHL] pour la d\émonstration des autre points, et quelques exemples explicites. 
 
\subsection{Exemples pour   $n=4$ :} 

Dans ce cas,  
 \n  $d_1=\frac{   1}{6}(k_0+2)(k_0+3)(k_0+4) \ \bigl(=c(4,k_0+1)\bigr)$, \hb 
 $d_2=\frac{ 1 }{6} (k_0+1)(k_0+3)(k_0+4) $, \hb  $d_3=\frac{  1}{6}(k_0+1)(k_0+2)(k_0+4)$,\hb  $d_4=\frac{   1}{6}(k_0+1)(k_0+2)(k_0+3) \ \bigl(=c(4,k_0)\bigr)$.
 
 \n D\éfinissant  $P=\partial _0$ en prenant  les mineurs de taille maximale de $M$,      la suite
 %$0\to \oplus _{k=1}^K{\cal O}_{\hbox {\bb P} }(-c_k) \buildrel N \over \longrightarrow \oplus _{j=1}^J{\cal O}_{\hbox {\bb P} }(-b_j) \buildrel M \over \longrightarrow \oplus _{i=1}^I{\cal O}_{\hbox {\bb P} }(-a_i) $
  $$0\to E_3 \buildrel N \over \longrightarrow E_2 \buildrel M \over \longrightarrow E_1\buildrel P \over \longrightarrow E_0\cong {\cal O}_{\hbox {\bb P} }(0)$$ s'\écrit :
$$0\to (d-d_4){\cal O}_{\hbox {\bb P} }(-(k_0+4))  \buildrel N \over \longrightarrow   3 (d-d_3){\cal O}_{\hbox {\bb P} }(-(k_0+3)) \buildrel M \over \longrightarrow   3 (d-d_2){\cal O}_{\hbox {\bb P} }(-(k_0+2)) \oplus\ (d_1-d){\cal O}_{\hbox {\bb P} }(-(k_0+1))\buildrel P \over \longrightarrow  {\cal O}_{\hbox {\bb P} }(0)  ,$$ 
$$0\to (d-d_4){\cal O}_{\hbox {\bb P} }(-(k_0+4))  \buildrel N \over \longrightarrow   3(d-d_3){\cal O}_{\hbox {\bb P} }(-(k_0+3))\oplus\ 3(d_2-d){\cal O}_{\hbox {\bb P} }(-(k_0+2)) \buildrel M \over \longrightarrow (d_1-d){\cal O}_{\hbox {\bb P} }(-(k_0+1))\buildrel P \over \longrightarrow  {\cal O}_{\hbox {\bb P} }(0),$$
 $$0\to (d-d_4){\cal O}_{\hbox {\bb P} }(-(k_0+4))  \buildrel N \over \longrightarrow   3(d-d_3){\cal O}_{\hbox {\bb P} }(-(k_0+3))\oplus\ 3(d_2-d){\cal O}_{\hbox {\bb P} }(-(k_0+2)) \buildrel M \over \longrightarrow (d_1-d){\cal O}_{\hbox {\bb P} }(-(k_0+1))\buildrel P \over \longrightarrow  {\cal O}_{\hbox {\bb P} }(0), $$ 
selon que  $d\in [d_2,d_1[ $, $d\in [d_3,d_2] $, ou $d\in [d_4,d_3]. $

 \n Notons $(X,Y,Z,T,U)$ les coordonn\ées homog\ènes dans {\bb P}$_4$, et $(t,s)$ dans {\bb P}$_1$.
    
  \n   1) Pour  $k_0=1$ et $d=d_3$,  la courbe  monomiale   $$(t,s)\mapsto (X=t^2s^3,\ Y=t^3s^2,\ Z=t^4s,\ T=t^5,\ U=s^5)$$     est ordinaire et  acm. Son  id\éal $\cal I$  peut \^etre d\éfini par cinq   g\én\érateurs de  degr\é deux, et la  r\ésolution  suivante :
  $$0\to  {\cal O}_{\hbox {\bb P} }(-5) \buildrel N \over \longrightarrow  {\cal O}^{\oplus 5}   _{\hbox {\bb P} }(-3) \buildrel M \over \longrightarrow  {\cal O}^{\oplus 5}   _{\hbox {\bb P} }(-2)\buildrel P \over \longrightarrow  {\cal O}_{\hbox {\bb P} }\to {\cal O}_\Gamma\to 0$$
  avec
  $$P=(\!(X^2-ZU,XY-TU,Y^2-XZ,YZ-XT,Z^2-YT)\!)\ ,M =\begin{pmatrix}
  0&0&-T&Z&-Y\\
  0&0&Z&-Y&X\\
  T&-Z&0&X&0\\
  -Z&Y&-X&0&-U\\
  Y&-X&0&U&0\\
  \end{pmatrix}, \ N=P^*.$$

\n 2)  Pour  $k_0=2$ et $d=d_2$, la courbe  monomiale   $$(t,s)\mapsto (X=t^8s^7,\ Y=t^{12}s^3,\ Z=t^{13}s^2,\ T=t^{15},\ U=s^{15})$$     est ordinaire et  acm. 
Son  id\éal $\cal I$  peut \^etre d\éfini par cinq   g\én\érateurs de  degr\é trois,   deux de  degr\é  quatre, et la  r\ésolution  suivante  (minimale, bien que les in\égalit\és  entre les $a_{j,i}$ ne soient pas strictes, la sous-matrice $2\times 2$ de $M$ dans le coin en haut \à droite ne contenant que des z\éros) : 
 $$0\to 5\ {\cal O}_{\hbox {\bb P} }(-6) \buildrel N \over \longrightarrow 9\ {\cal O}_{\hbox {\bb P} }(-5)\oplus 2\ {\cal O}_{\hbox {\bb P} }(-4) \buildrel M \over \longrightarrow 2\ {\cal O}_{\hbox {\bb P} }(-4) \oplus 5\ {\cal O}_{\hbox {\bb P} }(-3)\buildrel P \over \longrightarrow {\cal O}_{\hbox {\bb P} }\to {\cal O}_{ \Gamma}\to 0,    $$
 avec 
$$  P=(\!( \   X^2Z^2-YT^2U\  ,\ XY^2Z-T^3U\ ,\ X^3-Y^2U\ ,\  X^2Y-ZTU\ ,\   Y^3-XZT\ ,\   YZ^2-XT^2\ ,\   Z^3-Y^2T\   )\!), $$ 
$$M=\begin{pmatrix}
0&0&0&Z&0&T&Y&0&-X& 0\hskip .5cm \ 0 \\

0&0&0&0&X&-Y&0&Z&0& 0\hskip .5cm \ 0\\

-ZT&0&T^2&0&0&0&0&0&Z^2&Y \hskip .5cm  0\\

Y^2&0&-Z^2&-YT&-YZ&0&-Z^2&-T^2&0& \!\!\!\!\!-X\hskip .5cm  0\\

-X^2&Z^2&0&0&0&XZ&0&0&0&U \hskip .2cm -T \\

0&-Y^2&X^2&0&-TU&0&0&-XY&YU&0 \hskip .5cm Z\\

0&XT&-TU&-X^2&0&0&-TU&0&0&0 \hskip .2cm -Y\\
 \end{pmatrix}$$
%\hbox
  {et } 
 $$ N=\begin{pmatrix}

 Z&0&T&0&0\cr
 0 &X&0&-U&0\cr
 0  &0&Z&0&Y\cr
  0  &T&Y&0&0\cr
   Y  &0&0&0&-Z\cr
    X  &-Z&0&0&0\cr
      0 &0&-Z&-X&0\cr
       0 &-Y&0&0&X\cr
        T &0&0&-Y&0\cr
         0 &0&0&Z^2&-T^2\cr
           0&0&-X^2&0&-TU\cr
           
    \end{pmatrix}.$$
    
\n {\bf R\éf\érences :}

\noindent [CL] V. Cavalier et D. Lehmann, Ordinary holomorphic webs of codimension one, preprint  (arXiv  math 0703596 v2[mathDS], 13/10/2008),   Ann. Sc. Norm. Super.  Pisa, cl. Sci (5), vol XI (2012), 197-214.

 %\pagebreak

%\noindent [CL2]   V. Cavalier et D. Lehmann, Rang et courbure de Blaschke des tissus holomorphes r\éguliers de codimension un, C.R. Acad. Sci, Paris, Ser.I 346 (2008).

\noindent [E] D. Eisenbud, Commutative Algebra with a view Toward Algebraic Geometry, Graduate texts in Mathematics 150, Springer, 1994)

\noindent [El] G. Ellingsud, Sur le sch\éma de Hilbert des vari\ét\és de codimension 2 dans {\bb P}$_e$ \à c\^one de Cohen-Macaulay, Ann. Sc. Ec. Norm. Sup., t. 8, fasc. 4, 423-431, 1975)

\noindent [GHL] L. Gruson, Y. Hantout  et  D. Lehmann, Courbes alg\ébriques ordinaires et tissus associ\és, C.R.Acad.Sci.Paris, Ser.I,  350 (2012), 513-518.

\noindent [GP] L. Gruson, et  C. Peskine, Genre des courbes de l'espace projectif, Algebraic Geometry, Tromso$\!\!\!/$ 1977, Springer Verlag, Lecture Notes in Math.   687, 31-59, 1978.

\noindent [H]   J. Harris,  Curves in projective space, chapter III (with the collaboration of D. Eisenbud), Les Presses de l'Universit\é de Montr\éal, 1982.

\bigskip
  % \n Laurent Gruson,  {Maths., Universit\é   de Versailles, 45 Avenue des Etats Unis, 78035 Versailles, France} \hb  email : {laurent.gruson@math.uvsq.fr}, 
    \n Youssef Hantout,  {Dept.  Math\ématiques, Universit\é de Lille 1, 59650 Villeneuve d'Ascq cedex, France}\hb email : {hantout@math.univ-lille1.fr}, 
   \n Daniel Lehmann,  {4 rue Becagrun,  30980 Saint Dionisy, France}\hb  email : {lehm.dan@gmail.com},

 \end{document}